\documentclass{amsart}
\usepackage{graphicx,amsfonts,amssymb,amsmath,amsthm}
\usepackage[pdftex]{hyperref}
\usepackage{cite}

\theoremstyle{plain} 
\newtheorem{theorem}    {Theorem}[section]

\newtheorem{question}   [theorem]{Question}

\theoremstyle{definition}

\theoremstyle{remark}
\newtheorem{remark}              {Remark}

\numberwithin{equation}{section}

\def\C{\mathbb C}

\def\Q{\mathbb Q}
\def\R{\mathbb R}

\begin{document}

\title[Distribution of Hecke eigenvalues for {\rm GL}(2)]{On the distribution of Hecke eigenvalues for cuspidal automorphic representations for {\rm GL}(2)}

\author{Nahid Walji}
\address{Institut f\"ur Mathematik, Universit\"at Z\"urich, Winterthurerstrasse 190, CH-8057, Z\"urich, Switzerland}
\email{nahid.walji@math.uzh.ch}
\maketitle
\begin{abstract} Given a self-dual cuspidal automorphic representation for {\rm GL}(2) over a number field, we establish the existence of an infinite number of Hecke eigenvalues that are greater than an explicit positive constant, and an infinite number of Hecke eigenvalues that are less than an explicit negative constant. This provides an answer to a question of Serre. We also consider analogous problems for cuspidal automorphic representations that are not self-dual. \end{abstract}

\section{Introduction}

Let $\pi$ be a cuspidal automorphic representation for {\rm GL}(2)/$F$, where $F$ is a number field. Associated to a finite place $v$ at which $\pi$ is not ramified is the Hecke eigenvalue $a_v(\pi)$ of $\pi$. Given the sequence $\{a_v(\pi)\}$ of Hecke eigenvalues of $\pi$, one can ask how these are distributed.
\begin{question}[Serre]\label{q1}
Let $\pi$ be a self-dual cuspidal automorphic representation for {\rm GL}(2)/ $\Q$. Can we find positive constants $c, c'$ such that, for any $\epsilon > 0$, we can show the existence of infinitely many primes $p$ such that $a_p(\pi) > c - \epsilon$ as well as infinitely many primes $p$ such that $a_p(\pi) < c' + \epsilon$?
\end{question}

This question was posed by Serre in his appendix to a paper of Shahidi ~\cite{Sh94}. In this appendix, he determines such one-sided bounds in the case of holomorphic cusp forms with constants $c$ and $c'$ both equal to $2\cos (2 \pi /7)$. He then determines upper bounds (of the same size) on the absolute value of Hecke eigenvalues associated to Maass forms, and finishes that subsection by raising the question of whether it is possible to obtain one-sided bounds (of the form as those in Question ~\ref{q1}) in the case of Maass forms. He points out in his method that an obstruction is the lack of knowledge of the Ramanujan conjecture (or, at least, bounded Hecke eigenvalues) for Maass forms.

Whilst assuming the Ramanujan conjecture (or at least that the Satake parameters are bounded), Ramakrishnan~\cite{Ra97} obtained related results on the distribution of Hecke eigenvalues for cuspidal automorphic representations for {\rm GL}(n) over a number field.
There is also the older work of Murty where he shows the existence of infinitely many large absolute values of Hecke eigenvalues in the setting of holomorphic modular forms ~\cite{Mu83}. 
More recently in 2002, Kim and Shahidi applied their cuspidality results ~\cite{KS02} to increase Serre's bounds (using the same approach as Serre) to $2 \cos (2 \pi /11)$.
When $F$ is a totally real field, one now knows more about the distribution of Hecke eigenvalues for (non-CM) holomorphic Hilbert modular forms due to the proof of the Sato-Tate conjecture by Barnet-Lamb, Gee, and Geraghty ~\cite{BLGG11}.

The case of Maass forms, on the other hand, does not appear to be accessible through geometric methods, and our aim is to make progress in this case:

\begin{theorem} \label{t1}
Let $\pi$ be a self-dual cuspidal automorphic representation for ${\rm GL}(2)$ over a number field. Then:\\
(i) For any $\epsilon > 0$, there exist infinitely many places $v$ such that 
\begin{align*}
a_v (\pi) > 0.904... - \epsilon.
\end{align*}
(ii) If $\pi$ is not monomial, then for any $\epsilon > 0$, there exist infinitely many places $v$ such that 
\begin{align*}
a_v (\pi) < -1.164... + \epsilon.
\end{align*}
\end{theorem}

\begin{remark}
(a) In part (ii) of the theorem above, the non-monomial condition is necessary. For example, consider the unique 2-dimensional irreducible representation of $S_3$. The character of this representation takes the values 2, 0, and $-1$, and therefore the inequality in part (ii) would not hold.\\
(b) The bound in part (ii) has larger absolute value than the bound in part (i) because we are able to exploit positivity and so make use of higher symmetric power $L$-functions.
\end{remark}

The proof of the theorem above makes use of results on the automorphy of symmetric powers, as well as information on the analytic properties of the $L$-functions of higher symmetric powers.
These include the automorphy and cuspidality results of Gelbart-Jacquet ~\cite{GJ78}, Kim-Shahidi ~\cite{KS00,KS02}, and Kim ~\cite{Ki03}. 
Our bounds are not merely a modern update of Serre's method; for example, even if using results only available at the time of ~\cite{Sh94}, our approach still produces non-trivial one-sided bounds for Maass forms.

In the case of cuspidal automorphic representations that are not self-dual, the Hecke eigenvalues will not all be real (and so we cannot ask for bounds as in the form of Question ~\ref{q1}). Instead we take the real part of the Hecke eigenvalue and seek to determine the existence of an infinite number of Hecke eigenvalues that lie a certain distance to the left, or to the right, of the imaginary axis. We also generalise this by replacing the imaginary axis with any complex line, denoted as $e ^{i (\pi/2-\phi)}\R$, through the origin.

\begin{theorem}\label{t2}
Let $\pi$ be a non-self-dual unitary cuspidal automorphic representation for ${\rm GL}(2)$ over a number field, that is non-monomial. Fix any $\phi \in [0,\pi]$.
Then for any $\epsilon > 0$, there exist infinitely many places $v$ such that 
\begin{align*}
{\rm Re}(a_v(\pi)  e ^{i\phi}) > 0.5 - \epsilon,
\end{align*}
and infinitely many places $v$ such that 
\begin{align*}
{\rm Re}(a_v(\pi) e ^{i\phi}) < -0.5 + \epsilon.
\end{align*}
\end{theorem}

The theorems above are also related to the work of Kohnen, Knopp, and Pribitkin on sign changes for holomorphic modular forms ~\cite{KKP03}. These authors show the existence of an infinite number of Fourier coefficients lying on each side of the imaginary axis in $\C$, and also prove this for all other complex lines through the origin. A consequence of Theorem ~\ref{t2} is a generalisation of this result to Maass forms, and over any number field.

The rest of the paper is structured as follows: Section ~\ref{spre} recalls some necessary background as well as proving technical results that will be needed later. Section ~\ref{ssd} consists of the proof for the self-dual case and Section ~\ref{snsd} consists of the proof for the non-self-dual case. 

We will prove the theorems for when $F = \Q$. The general case follows in the same way.

\section{Preliminaries} \label{spre}

\subsection{Lower and upper Dirichlet density}

Recall that the lower Dirichlet density for a set $S$ of primes is defined to be 
\begin{align*}
\underline{\delta}(S) := \lim_{s \rightarrow 1 ^+}\inf \frac{\sum_{p \in S}p ^{-s}}{-\log (1/ (s-1))}
\end{align*}
and the upper Dirichlet density is 
\begin{align*}
\overline{\delta}(S) := \lim_{s \rightarrow 1 ^+}\sup \frac{\sum_{p \in S}p ^{-s}}{-\log (1/ (s-1))}.
\end{align*}

When they are equal, then the set $S$ is said to possess a Dirichlet density, which is equal to the upper/lower Dirichlet density.\\

We also make note of some identities involving the limit supremum and infimum that will be of use later on.
Given real-valued functions $A(s), B(s)$ and a point $t \in \R$, we have 
\begin{align*}
\lim_{s \rightarrow t}{\rm sup} (A(s) + B(s)) &\leq \lim_{s \rightarrow t}{\rm sup}\ A(s) + \lim_{s \rightarrow t}{\rm sup}\ B(s) \\
 - \lim_{s \rightarrow t}{\rm inf}(-A(s)) &= \lim_{s \rightarrow t}{\rm sup}\ A(s).
\end{align*}
If $\lim_{s \rightarrow t} (A(s) + B(s))$ exists, then 
\begin{align}\label{infsup}
\lim_{s \rightarrow t}{\rm inf}(A(s)) + \lim_{s \rightarrow t}{\rm sup}(B(s)) = \lim_{s \rightarrow t} (A(s) + B(s)).
\end{align}
Lastly, if $A(s), B(s) \geq 0$, then 
\begin{align*}
\lim_{s \rightarrow t}{\rm sup}(A(s) \cdot B(s)) &\leq \lim_{s \rightarrow t}{\rm sup}(A(s)) \cdot \lim_{s \rightarrow t}{\rm sup}(B(s)).\\
\end{align*}

\subsection{Asymptotic behaviour of Dirichlet series}\label{ds}

Assume that $\pi$ is non-monomial. We want to understand the asymptotic behaviour of the Dirichlet series $\sum_{p}\frac{a_p(\pi) ^k}{p^s}$ for $k = 2,3,4$, as $s \rightarrow 1^+$. We will also address what can be said in the case of $k = 6,7$ and $8$.

When $k = 2$, one knows that the Rankin--Selberg $L$-function $L(s,\pi \times \pi)$ is holomorphic in $\C$, except for potentially a simple pole at $s=1$ iff $\pi$ is self-dual. We make use of bounds towards the Ramanujan conjecture, by Kim--Sarnak in the rational case (see Appendix 2 of~\cite{Ki03}) and Blomer--Brumley ~\cite{BB11} in the number field case. Thus we know that the size of the Satake parameters at each prime $p$ is bounded above by $p ^{7/64}$, which implies that 
\begin{align*}
\log L(s,\pi \times \pi)= \sum_{p}\frac{a_p(\pi) ^2}{p^s} + O\left(1\right)
\end{align*}
as $s \rightarrow 1^+$, and thus either 
\begin{align*}
\sum_{p}\frac{a_p(\pi) ^2}{p^s} \sim \log \left(\frac{1}{s-1}\right)
\end{align*}
when $\pi$ is self-dual, or 
\begin{align*}
\sum_{p}\frac{a_p(\pi) ^2}{p^s} = O\left(1\right)
\end{align*}
when $\pi$ is not self-dual.\\

For $k = 3$, we observe that, for a set $T$ containing exactly all the primes at which $\pi$ is ramified and the archimedean place, 
\begin{align*}
L^T(s, \pi \times \pi \times \pi) = L^T(s, \pi, {\rm Sym}^3) L^T(s, \pi \otimes \omega)^2
\end{align*}
due to the Klebsch--Gordon decomposition of tensor powers, where $\omega$ is the central character of $\pi$. From Shahidi~\cite{Sh94}, we know that $L(s,\pi, {\rm Sym}^3)$
is invertible at $s=1$, and so the same holds for $L^T(s, \pi, {\rm Sym}^3)$ given the invertibility of the finite number of factors associated to places in $T$.
The same holds for $L^T(s, \pi \otimes \omega)^2$ and thus for $L^T(s, \pi \times \pi \times \pi)$. Bounds towards the Ramanujan conjecture then imply that (whether or not $\pi$ is self-dual) 
\begin{align*}
\sum_{p}\frac{a_p(\pi) ^3}{p^s} = \log L^T(s, \pi \times \pi \times \pi) + O\left(1\right) = O\left(1\right).\\
\end{align*}

For $k = 4$, we note that 
\begin{align*}
L^T(s, \pi ^{\times 4}) = L^T(s, \pi, {\rm Sym}^4) L^T(s, \pi, {\rm Sym}^2 \otimes \omega)^3 L^T(s, \omega ^2)^2.
\end{align*}
Now, $L^T(s, \pi, {\rm Sym}^4)$ is invertible at $s=1$~\cite{Sh94} and so is $L^T(s, \pi, {\rm Sym}^2 \otimes \omega)$. We only need the case where $\pi$ is self-dual, and so the central character $\omega$  is trivial, implying that $L(s,\pi ^{\times 4})$ has a pole of order two at $s=1$. Bounds towards the Ramanujan conjecture then imply 
\begin{align} \label{k4}
\sum_{p}\frac{a_p(\pi) ^4}{p^s} = \log L^T(s, \pi ^{\times 4}) + O\left(1\right)
\sim 2 \log \left(\frac{1}{s-1}\right).\\ \notag 
\end{align}

\noindent \textbf{Case $k = 6$:} 
Assume that $\pi$ is not of dihedral or tetrahedral type. Then we know that ${\rm Sym}^3 \pi$ is cuspidal, due to Kim and Shahidi ~\cite{KS00,KS02}. We note 
\begin{align*}
L^T(s, \pi ^{\times 6}) = L^T(s, {\rm Sym}^3 \pi \times {\rm Sym}^3 \pi) 
L^T(s, {\rm Sym}^3 \pi \times \pi \otimes \omega )^4 
L^T(s, \pi \times \pi \otimes \omega ^2)^4.
\end{align*}
Therefore, when $\pi$ is self-dual, $L^T(s, \pi ^{\times 6})$ has a pole of order 5 at $s=1$.

If $\pi$ is of tetrahedral type, then ${\rm Sym}^3 \pi$ is not cuspidal. It has the following decomposition 
\begin{align}\label{k6d}
{\rm Sym}^3 \pi \otimes \omega ^{-1} \simeq (\pi \otimes \mu) \boxplus (\pi \otimes \mu ^2),
\end{align}
where $\mu$ is a non-trivial Hecke character of order three. Thus, for self-dual $\pi$, we have 
\begin{align*}
L^T(s, {\rm Sym}^3 \pi \times {\rm Sym}^3 \pi) = L^T(s, \pi \times \pi \otimes \mu ^2) L^T(s, \pi \times \pi \otimes \mu ^3)^2 L^T(s, \pi \times \pi \otimes \mu),
\end{align*}
and so we conclude that $L^T(s, \pi ^{\times 6})$ has a pole of order 6 at $s=1$.

Now write $\{\alpha_p, \beta_p\}$ for the Satake parameters of $\pi$ at $p$. For non-dihedral self-dual $\pi$, bounds towards the Ramanujan conjecture imply that 
\begin{align}
\label{s6} \sum_{k = 1 }^{ 4} \sum_{p} \frac{(\alpha_p ^k + \beta_p ^k)^6}{p^s} \geq 5 \log \left(\frac{1}{s-1}\right) + o\left(\log \left(\frac{1}{s-1}\right)\right).\\ \notag 
\end{align}

\noindent \textbf{Case $k = 7$:}
Assume that $\pi$ is not of solvable polyhedral type (that is, it is not dihedral, tetrahedral, or octahedral). Then we know that ${\rm Sym}^3 \pi$ and ${\rm Sym}^4 (\pi)$ are cuspidal automorphic representations ~\cite{KS02}. Now 
\begin{align}\label{k7}
L^T(s, \pi ^{\times 7}) = L^T(s, {\rm Sym}^4 \pi \times {\rm Sym}^3 \pi) L^T(s, {\rm Sym}^4 \pi \times \pi \otimes \omega) ^2 \\ \notag 
\cdot L^T(s, {\rm Sym}^3 \pi \times {\rm Sym}^2 \pi \otimes \omega)^3 L^T(s, {\rm Sym}^3 \pi \otimes \omega ^3)^2 \\ \notag 
\cdot L^T(s, {\rm Sym}^2 \pi \times \pi \otimes \omega ^2) ^6 L^T(s, \pi \otimes \omega ^3)^4.
\end{align}
Since all the $L$-functions on the right-hand side are invertible at $s=1$, the same holds for $L^T(s, \pi ^{\times 7})$.

If $\pi$ is octahedral, then ${\rm Sym}^3 \pi$ is still cuspidal, but ${\rm Sym}^4 \pi$ is not. Its decomposition is 
\begin{align*}
{\rm Sym}^4 \pi \otimes \omega ^{-1} \simeq (\pi (\chi ^{-1}) \otimes \omega) \boxplus {\rm Sym}^2 \pi,
\end{align*}
where $\pi (\chi ^{-1})$ is the monomial representation of {\rm GL}(2)/$\Q$ associated to a certain non-trivial Hecke character $\chi$ (see Section 3 of ~\cite{KS02}). Applying this decomposition to equation ~\ref{k7}, we again see that $L^T(s, \pi ^{\times 7})$ is invertible at $s=1$.
If $\pi$ is tetrahedral, ${\rm Sym}^4 \pi$ is not cuspidal, and decomposes as 
\begin{align*}
{\rm Sym}^4 \pi \otimes \omega ^{-1} \simeq {\rm Sym}^2 \pi \boxplus \omega \mu \boxplus \omega \mu ^2,
\end{align*}
where $\mu$ is the same Hecke character as in ~\ref{k6d}. Using the equation above as well as equation ~\ref{k6d} in equation ~\ref{k7} implies that $L^T(s, \pi ^{\times 7})$ is invertible at $s=1$.\\

\noindent \textbf{Case $k = 8$:}
If $\pi$ is not of solvable polyhedral type (that is, it is not dihedral, tetrahedral, or octahedral), then Kim ~\cite{Ki03} has shown that ${\rm Sym}^4 (\pi)$ is a cuspidal automorphic representation. Therefore, the associated $L$-function $L(s, {\rm Sym}^4 (\pi))$ is invertible at $s=1$.

Let $T$ be the set containing exactly the primes at which $\pi$ is ramified, as well as the archimedean place. Clebsch-Gordon decompositions show that 
\begin{align*}
L^T(s, \pi ^{\times 8}) =& L^T(s, {\rm Sym}^4 \pi \times {\rm Sym}^4 \pi) 
L^T(s, {\rm Sym}^4 \pi \times {\rm Sym}^2\pi \otimes \omega)^6 \\
\cdot & 
L^T(s, {\rm Sym}^2 \pi \otimes \omega \times {\rm Sym}^2\pi \otimes \omega)^9
L^T(s, {\rm Sym}^4 \pi \otimes \omega ^2)^4 \\
\cdot &
L^T(s, {\rm Sym}^2 \pi \otimes \omega ^3)^{12}
L^T(s, \omega ^4)^4,
\end{align*}
where $\omega$ is the central character of $\pi$. We conclude that $L^T(s, \pi ^{\times 8})$ has a pole of order 14 at $s=1$. 
Using the known bounds towards the Ramanujan conjecture, we get
\begin{align}\label{s8}
\sum_{m = 1 }^{ 8}\sum_{p \not \in T}\frac{(\alpha_p^m + \beta_p ^m)^8}{p^{ms}} = 14 \log \left(\frac{1}{s-1}\right) + o\left(\log \left(\frac{1}{s-1}\right)\right).\\ \notag
\end{align}

\begin{remark}
One might ask why we do not also make use of $L(s,\pi, {\rm Sym}^9)$, given that we know (by Kim and Shahidi ~\cite{KS06}) that for a cuspidal automorphic representation for {\rm GL}(2) over a number field that is not of solvable polyhedral type, the $L$-function $L(s, \pi, {\rm Sym}^9)$ is holomorphic in the interval $(1,2]$, and at $s=1$ the $L$-function either has a simple pole, is invertible, or has a simple zero. The issue is that we do not know of the absolute convergence of the Euler product for ${\rm Re}(s)> 1$ (it is only known for ${\rm Re}(s)> 2$), which is required by this method.
\end{remark}

\subsection{Positive coefficients in Dirichlet series:}\ 

Now we note that, for self-dual $\pi$ (and so the central character is trivial)
\begin{align*}
a_p ^m = (\alpha_p ^m + \beta_p ^m) + ^m C_1  (\alpha_p ^{m-1} + \beta_p ^{m-1}) +     \dots 
+ ^m C_{m/2 + 1}(\alpha_p ^2 + \beta_p ^2) + ^m C_{m/2}\cdot 1,
\end{align*}
for even $m$, and 
\begin{align*}
a_p ^m = (\alpha_p ^m + \beta_p ^m) + ^m C_1  (\alpha_p ^{m-1} + \beta_p ^{m-1}) +     \dots 
+ ^m C_{\frac{m + 3}{2}}(\alpha_p ^3 + \beta_p ^3) + ^m C_{\frac{m + 1}{2}} (\alpha_p + \beta_p ),
\end{align*}
for odd $m$. Then $(\alpha_p ^m + \beta_p ^m)$ can be written as a linear combination of $a_p ^m , \dots, a_p ^0$. Since $a_p$ is real, so is $(\alpha_p ^m + \beta_p ^m)$ and thus $(\alpha_p ^m + \beta_p ^m)^6$ and $(\alpha_p ^m + \beta_p ^m)^8$ are non-negative for all $p$ and $m$. In particular,
\begin{align}\label{s8ii}
\sum_{p}\frac{|a_p|^8}{p^s} \leq 
14 \log \left(\frac{1}{s-1}\right) + o\left(\log \left(\frac{1}{s-1}\right)\right).
\end{align}

In the case of $(\alpha_p ^m + \beta_p ^m)^7$, we take a different approach. Note that if $a_p$ is sufficiently large (say $a_p > c$ for some constant $c$), then $(\alpha_p ^m + \beta_p ^m) > 0$ for $m = 1, \dots ,4$. We define 
\begin{align}\label{def1}\notag 
A &:= \{p \mid a_p > 0\}\\ \notag
A_b &:= \{p \in A \mid a_p \leq c\}\\
A_u &:= \{p \in A \mid a_p > c\}.
\end{align}
Then we note that for
\begin{align*}
\sum_{m = 1 }^{ 4} \sum_{p \in A} \frac{(\alpha_p ^m + \beta_p ^m)^7}{p ^{ms}}
= \sum_{p \in A}\frac{a_p ^7}{p^s} + \sum_{m = 2 }^{ 4} \sum_{p \in A_u}\frac{(\alpha_p ^m + \beta_p ^m)^7}{p ^{ms}} +  \sum_{m = 2 }^{4} \sum_{p \in A_b}\frac{(\alpha_p ^m + \beta_p ^m)^7}{p ^{ms}},
\end{align*}
the last double sum on the right-hand side is bounded as $s \rightarrow 1^+$, and the other two series on the right-hand side all have positive coefficients.

\section{Large positive and negative eigenvalues in self-dual case} \label{ssd}

\subsection{Positive side}\label{p1pos}\ 
If $\pi$ is a cuspidal automorphic representation for {\rm GL}(2)/$\Q$ that is of solvable polyhedral type, then it is associated to an Artin representation, in which case we immediately know that the bound in part (i) of Theorem ~\ref{t1} holds.
So let us assume that $\pi$ is not of solvable polyhedral type. 

Let $A$ consist of the set of primes $p$ for which the associated Hecke eigenvalues are positive, and $B$ will consist of the set of primes $p$ for which $a_p \leq 0$.

From equation ~\ref{k4} we know that   
\begin{align*}
\lim_{s \rightarrow 1 ^+}{\rm sup} \frac{\sum_{p}\frac{|a_p|^4}{p^s}}{\log \left(\frac{1}{s-1}\right)} = 2.
\end{align*}
Therefore, 
\begin{align} \label{eqn1}
\lim_{s \rightarrow 1 ^+}{\rm sup} \frac{\sum_{p \in A}\frac{|a_p|^4}{p^s}}{\log \left(\frac{1}{s-1}\right)} + \lim_{s \rightarrow 1 ^+}{\rm sup} \frac{\sum_{p \in B}\frac{|a_p|^4}{p^s}}{\log \left(\frac{1}{s-1}\right)} \geq 2.
\end{align}
Let us denote
\begin{align*}
d := \lim_{s \rightarrow 1 ^+}{\rm sup} \frac{\sum_{p \in B}\frac{|a_p|^4}{p^s}}{\log \left(\frac{1}{s-1}\right)}
\end{align*}
Now 
\begin{align*}
\sum_{p \in B}\frac{|a_p|^{8/5}|a_p|^{12/5}}{p^s} \leq \left(\sum_{p \in B}\frac{|a_p|^8}{p^s}\right)^{1/5}
\left(\sum_{p \in B}\frac{|a_p|^3}{p^s}\right)^{4/5}.
\end{align*}

We have, by equation ~\ref{s8ii},
\begin{align*}
\lim_{s \rightarrow 1 ^+}{\rm sup} \frac{\sum_{p \in B}\frac{|a_p|^8}{p^s}}{\log \left(\frac{1}{s-1}\right)} \leq 14.
\end{align*}
So dividing the inequality above by $\log (1 / (s-1))$ and taking the limit supremum as $s \rightarrow 1^+$
\begin{align*}
d &\leq 14 ^{1/5} \left(\lim_{} {\rm sup} \frac{\sum_{p \in B}\frac{|a_p|^3}{p^s}}{\ell (s)}\right)^{4/5}\\
\frac{d ^{5/4}}{14 ^{1/4}} &\leq \lim_{}{\rm sup} \frac{\sum_{p \in B}\frac{|a_p|^3}{p^s}}{\log \left(\frac{1}{s-1}\right)}.
\end{align*}
Now 
\begin{align*}
\sum_{p}\frac{a_p ^3}{p^s} = 
 O\left(1\right)
\end{align*}
as $s \rightarrow 1 ^+$, and since the series converges absolutely, we write 
\begin{align*}
\sum_{p \in A} \frac{a_p ^3}{p^s} = - \sum_{p \in B} \frac{a_p ^3}{p^s} + O\left(1\right).
\end{align*}
and so 
\begin{align} \label{s1}
\lim_{s \rightarrow 1 ^+}{\rm sup} \frac{\sum_{p \in A}\frac{|a_p|^3}{p^s}}{\log \left(\frac{1}{s-1}\right)} = \lim_{s \rightarrow 1 ^+}{\rm sup} \frac{\sum_{p \in B}\frac{|a_p|^3}{p^s}}{\log \left(\frac{1}{s-1}\right)} \geq \frac{d ^{5/4}}{14 ^{1/4}}.
\end{align}

Note that equation~\ref{eqn1} implies that 
\begin{align}\label{s2}
\lim_{s \rightarrow 1 ^+}{\rm sup} \frac{\sum_{p \in A}\frac{|a_p|^4}{p^s}}{\log \left(\frac{1}{s-1}\right)} \geq 2-d.
\end{align}
For any $\epsilon > 0$, equation ~\ref{s1} implies that there exists infinitely many $p \in A$ such that 
\begin{align*}
a_p > \frac{d ^{5/12}}{14 ^{1/12}} - \epsilon,
\end{align*}
and equation ~\ref{s2} implies that there exists infinitely many $p \in A$ such that 
\begin{align*}
a_p > \sqrt[4]{2-d} - \epsilon.
\end{align*}
We calculate 
\begin{align*}
{\rm min}_{d \in [0,2]} \left({\rm max}\left\{\frac{d ^{5/12}}{14 ^{1/12}}, \sqrt[4]{2-d}\right\}\right) = 0.904...,
\end{align*}
which occurs when $d = 1.331...$ .

We conclude that, for any positive $\epsilon$, there exists infinitely many $p$ such that 
\begin{align*}
a_p > 0.904... - \epsilon.
\end{align*}

\subsection{Negative side} \label{p1neg}\ 

Here, one could use the same method as in the previous subsection, but we can obtain a stronger bound by altering our approach and also making use of positivity.
(Note that the method in this section can be used to address the case of subsection ~\ref{p1pos}, but would result in a weaker bound of $1/ \sqrt{2} = 0.707...$, instead of $0.904...$.)

Assume that $\pi$ is a non-dihedral cuspidal automorphic representation for {\rm GL}(2)/$\Q$. Denote the Satake parameters at $p$ by $\{\alpha_p, \beta_p\}$.

Define set $A$ to consist of exactly those $p$ at which $a_p > 0$; and set $B$ to consist of exactly those $p$ at which $a_p \leq 0$.
Let us assume that (almost all) the Hecke eigenvalues associated to elements of set $B$ are bounded in size by some fixed $t > 0$.

First we establish two inequalities.\\

(i) Note that
\begin{align*}
\left|\sum_{k = 1}^{4} \sum_{B}\frac{(\alpha_p ^k + \beta_p ^k)^7}{p^{ks}}\right| 
\leq \sum_{B}\frac{t ^7}{p^s} + O\left(1\right)
\end{align*}
and 
\begin{align*}
\lim_{s \rightarrow 1^+ }{\rm sup} \left(\frac{\sum_{B}{t ^7}/{p^s}}{\ell (s)}\right) = \overline{\delta}(B) t ^7.
\end{align*}

The discussion in subsection ~\ref{ds} implies that
\begin{align*}
\sum_{k = 1}^{4}\sum_{p}\frac{(\alpha_p ^k + \beta_p ^k)^7}{p ^{ks}} =  O\left(1\right),
\end{align*}
and given the absolute convergence of the series for ${\rm Re}(s)> 1$, we write 
\begin{align*}
\sum_{k = 1}^{4}\sum_{p \in A}\frac{(\alpha_p ^k + \beta_p ^k)^7}{p ^{ks}} = 
- \sum_{k = 1}^{4}\sum_{p \in B}\frac{(\alpha_p ^k + \beta_p ^k)^7}{p ^{ks}} + O\left(1\right).
\end{align*}
So
\begin{align*}
\lim_{}{\rm sup} \left(\frac{\sum_{k = 1 }^{ 4}\sum_{p \in A}\frac{(\alpha_p ^k + \beta_p ^k)^7}{p ^{ks}}}{\ell (s)}\right)
&= \lim_{}{\rm sup} \left(\frac{\sum_{k = 1 }^{ 4}\sum_{p \in B}\frac{-(\alpha_p ^k + \beta_p^k)^7}{p ^{ks}}}{\ell (s)}\right)\\
&\leq \overline{\delta}(B) t ^7,
\end{align*}
and so, recalling the definitions in ~\ref{def1}, we have 
\begin{align}\label{bd1}
\lim_{}{\rm sup} \left( \frac{\sum_{p \in A} \frac{a_p ^7}{p^s} + 
\sum_{k = 1 }^{ 4}\sum_{p \in A_u}\frac{(\alpha_p ^k + \beta_p ^k)^7}{p ^{ks}}}{\ell (s)}\right)
&= \lim_{}{\rm sup} \left(\frac{\sum_{k = 1 }^{ 4}\sum_{p \in A}\frac{(\alpha_p ^k + \beta_p ^k)^7}{p ^{ks}}}{\ell (s)}\right)\\ \notag
&\leq \overline{\delta}(B) t ^7,
\end{align}
and note that all the coefficients in both series on the left-hand side are positive.\\

(ii) From equation ~\ref{s6} of subsection ~\ref{ds}, we have 
\begin{align*}
\lim_{}{\rm sup} \left(\frac{\sum_{k = 1 }^{ 4}\sum_{p}\frac{(\alpha_p ^k + \beta_p ^k)^6}{p ^{ks}}}{\ell (s)}\right) \geq 5,
\end{align*}
and so, 
\begin{align*}
\lim_{}{\rm sup} \left(\frac{\sum_{k = 1 }^{ 4}\sum_{p \in A}\frac{(\alpha_p ^k + \beta_p ^k)^6}{p ^{ks}}}{\ell (s)}\right) + \lim_{}{\rm sup} \left(\frac{\sum_{k = 1 }^{ 4}\sum_{p \in B}\frac{(\alpha_p ^k + \beta_p ^k)^6}{p ^{ks}}}{\ell (s)}\right) \geq 5,
\end{align*}
which gives 
\begin{align}\label{bd2}
\lim_{}{\rm sup} \left(\frac{\sum_{p \in A} \frac{a_p ^6}{p^s} + \sum_{k = 2}^{ 4}\sum_{p \in A_u}\frac{(\alpha_p ^k + \beta_p ^k)^6}{p ^{ks}}}{\ell (s)}\right)  \geq 5 - t ^6 \overline{\delta}(B).\\ \notag 
\end{align}

Now we combine the inequalities via the intermediary inequality (note that all the numerators of the fractions in the series are positive)
\begin{align*}
\left(\frac{\sum_{p \in A} \frac{a_p ^6}{p^s} + \sum_{k = 2 }^{ 4}\sum_{p \in A_u}\frac{(\alpha_p ^k + \beta_p ^k)^6}{p ^{ks}}}{\ell (s)}\right)
\leq &
\left(\frac{\sum_{p \in A} \frac{a_p ^7}{p^s}  + \sum_{k = 2 }^{ 4}\sum_{p \in A_u}\frac{(\alpha_p ^k + \beta_p ^k)^7}{p ^{ks}}}{\ell (s)}\right)^{6/7}\\
&\cdot \left(\frac{\sum_{p \in A} \frac{1}{p^s}  + \sum_{k = 2 }^{ 4}\sum_{p \in A_u}\frac{1}{p ^{ks}}}{\ell (s)}\right)^{1/7}.
\end{align*}
Taking the limit supremum of both sides and using the inequalities ~\ref{bd1},~\ref{bd2}, we obtain 
\begin{align*}
5 - t ^6 \overline{\delta}(B) &\leq \overline{\delta}(B)^{6/7} t ^6 \overline{\delta}(A)^{1/7}\\
5 - t ^6 &\leq t ^6
\end{align*}
which implies
\begin{align*}
1.164... &\leq t.
\end{align*}
We conclude that, for any $\epsilon > 0$, there exist infinitely many primes $p$ such that 
\begin{align*}
a_p < -1.164... + \epsilon.
\end{align*}

\begin{remark}
If the cuspidality of higher power symmetric lifts are known (along with stronger bounds towards the Ramanujan conjecture), then the lower bound on the size of $t$ arising from this method tends to 2, as expected.
\end{remark}

\section{Large positive and negative eigenvalues in non-self-dual case} \label{snsd}

Here we assume that $\pi$ is non-dihedral and not self-dual. Since the Hecke eigenvalues associated to $\pi$ will not all be real, we can consider the real parts of the Hecke eigenvalues instead. Rather than just studying the occurrence of Hecke eigenvalues that are of distance greater than some constant $c$ from one side or the other of the imaginary axis, we can choose any straight line in the complex plane that passes through the origin and study the occurrence of Hecke eigenvalues that are of distance at least $c$ from this line, first on one side, then the other.

Consider the two half-planes in $\C$ divided by the line $e ^{i(\pi/2 - \phi)}\R$. We determine the distance of a complex number $u$ from the line, and in which half-planes it lies, by examining the size and sign of ${\rm Re}(u e ^{i\phi })$.\\

At this point, the proof proceeds in a similar way to that in subsection ~\ref{p1neg}, though the slight variance in the method, along with the loss of self-duality, means that the bounds we obtain are different. In particular, we cannot use positivity to strengthen the bounds in this case (as we did in subsection ~\ref{p1neg}); as a result, note that the 5th to 9th symmetric power $L$-functions do not arise here.\\

We define the set $A$ to consist of exactly those primes $p$ such that ${\rm Re}(a_p e^{i\phi })> 0$. Set B will consist of exactly those primes $p$ where ${\rm Re}(a_p e^{i\phi})\leq 0$. Furthermore, we will assume that, for all but finitely many $p \in B$, $|{\rm Re}(a_p e^{i\phi})| \leq t$, for some $t$. We will write $\ell (s)$ to denote $\log (1 / (s-1))$.\\

We note that 
\begin{align*}
\lim_{} {\rm sup} \frac{\sum_{B} \frac{|{\rm Re}(a_p e^{i\phi})|^3}{p^s}}{\ell (s)} \leq \lim_{} {\rm sup} \frac{\sum_{B}\frac{t ^3}{p^s}}{\ell (s)} = t ^3 \overline{\delta}(B).
\end{align*}
Now 
\begin{align*}
\sum_{p} \frac{{\rm Re}(a_p e^{i\phi})^3}{p^s} &= \sum_{p}\frac{\left(\frac{1}{2}(a_p e^{i\phi}+ \overline{a_p}e^{-i\phi})\right)^3}{p^s} \\ &= 
\frac{1}{8} \left(\sum_{p }\frac{a_p ^3 e^{3i\phi} }{p^s} + 3 \sum_{p}\frac{a_p ^2 \overline{a_p}e^{i\phi}}{p^s} + 3 \sum_{p} \frac{a_p \overline{a_p}^2 e^{-i\phi} }{p^s} + \sum_{p} \frac{\overline{a_p}^3 e^{-3i\phi} }{p^s}\right)\\
&= o\left(\ell (s)\right),
\end{align*}
since each of these four series is $o\left(\ell (s)\right)$ as $s \rightarrow 1^+$.
So
\begin{align*}
\lim_{} {\rm sup} \frac{\sum_{A}  \frac{{\rm Re}(a_p e^{i\phi})^3}{p^s}}{\ell (s)}
+ \lim_{} {\rm inf} \frac{\sum_{B}  \frac{{\rm Re}(a_p e^{i\phi})^3}{p^s}}{\ell (s)} = 0 
\end{align*}
implying 
\begin{align} \label{nsd1}
\lim_{} {\rm sup} \frac{\sum_{A}  \frac{{\rm Re}(a_p e^{i\phi})^3}{p^s}}{\ell (s)}
 \leq  t ^3 \overline{\delta}(B).\\ \notag 
\end{align}

To obtain the second inequality, we start by observing
\begin{align*}
\sum_{p} \frac{{\rm Re}(a_p e^{i\phi})^2}{p^s} &=  \frac{1}{4} \left(\sum_{p }\frac{a_p ^2 e^{2i\phi} }{p^s} + 2 \sum_{p}\frac{a_p \overline{a_p}}{p^s} + \sum_{p} \frac{\overline{a_p}^2 e^{-2i\phi} }{p^s}\right)\\
&= \frac{1}{2}\ell (s) + o\left(\ell (s)\right),
\end{align*}
since the middle series contributes $2 \cdot (\ell (s) + o\left(\ell (s)\right))$, and the first and third series are $o\left( \ell (s) \right)$ since $\pi$ is not self-dual. And since
\begin{align*}
\lim_{} {\rm sup} \frac{\sum_{B}\frac{{\rm Re}(a_p e^{i\phi})^2}{p^s}}{\ell (s)} \leq 
\lim_{} {\rm sup} \frac{\sum_{B}\frac{t ^2}{p^s}}{\ell (s)} = t ^2 \overline{\delta}(B),
\end{align*}
we find that 
\begin{align}\label{nsd2}
\lim_{} {\rm sup} \frac{\sum_{A} \frac{{\rm Re}(a_p e^{i\phi})^2}{p^s}}{\ell (s)} \geq  \left(\frac{1}{2} - t ^2 \overline{\delta}(B)\right).
\end{align}

Now
\begin{align*}
\sum_{A} \frac{{\rm Re}(a_p e^{i\phi})^2}{p^s} \leq \left(\sum_{A} \frac{{\rm Re}(a_p e^{i\phi} )^3}{p^s}\right)^{2/3}\left(\sum_{A}\frac{1}{p^s}\right)^{1/3};
\end{align*}
dividing by $\ell (s)$, taking the limit supremum, and applying inequalities ~\ref{nsd1} and ~\ref{nsd2}
\begin{align*}
\frac{1}{2} - t ^2 \overline{\delta}(B) \leq 
(t ^3 \overline{\delta}(B)) ^{2/3} (\overline{\delta}(A))^{1/3}
\end{align*}
which implies that
\begin{align*}
1/2 \leq t.
\end{align*}
We conclude that, for any positive $\epsilon$, there exist infinitely many primes $p$ such that 
\begin{align*}
{\rm Re}(a_p e ^{i \phi}) < - 1/2 + \epsilon .\\
\end{align*}

To obtain the corresponding bound on the positive side, we again define sets $A$ and $B$ as before, but this time we assume that for almost all $p \in A$ we have $|{\rm Re}(a_p e ^{i\phi})| \leq t$. The proof then proceeds in the same manner as before.\\

\subsection*{Acknowledgements}
The author would like to thank Farrell Brumley, Emmanuel Kowalski, and Han Wu for their comments on an earlier version of this paper. 
This work began at UC Berkeley, and the author would like to thank Sug Woo Shin for hosting him as a Visiting Scholar.
At the University of Z\"urich, the author was supported by Forschungskredit grant K-71116-01-01 of the University of Z\"urich and partially supported by grant SNF PP00P2-138906 of the Swiss National Foundation.

\end{document}